\documentclass[12pt]{article}
\usepackage{amssymb}
\usepackage{amsmath,amsthm}

\providecommand{\U}[1]{\protect\rule{.1in}{.1in}}
\makeatletter
\def\@seccntformat#1{\csname the#1\endcsname.\quad}
\makeatother

\begin{document}

\title{{\Large \textbf{On a Variation of the Definition of Limit: Some Analytic Consequences}}}
\author{ Dhurjati Prasad Datta\thanks{
email:dp${_-}$datta@yahoo.com} \\
Department of Mathematics, University of North Bengal \\
Siliguri,West Bengal, Pin: 734013, India }
\date{}
\maketitle

\begin{abstract}
The basic formalism of a novel scale invarinat nonlinear analysis is presented. A few analytic number theoretic results are derived independent of standard approaches. 
\end{abstract}



\baselineskip=17.5pt
 
\section{Introduction}

Ordinary analysis is linear in the sense that a finite change (variation) in a real variable $x$ is formulated by reducing it to the level of linear elementary shifts (increments) of the form $x\mapsto X=x+h$, $h=X-x= \Delta x$ being the 'infinitesimal' elementary (that is to say, differential) increment. In this paper, we present basic building blocks of a {\em scale invariant, nonlinear}, analysis which would also support naturally (may be finite) elementary increments of the form $X/x\mapsto (x/X)^{\frac{1}{s}}, \ 0<s<1$. Since a study of Calculus must begin with a rigorous definition of limit, we also begin with a nontrivial definition of a {\em scale invariant limit}, by bringing in a {\em nonuniformity}  in the motion of a limiting variable, for instance, as $x$ approaches the limit 0. More precisely, we examine how a slight variation of conventional treatment and meaning of the statement $\underset{x \rightarrow 0} {\lim} \ x=0, \ x$ being a real variable can lead to potentially new results which would be significant for formulating a scale invariant analysis leading to a (dynamical) formation of Cantor sets profusely in the neighbourhood of $x=0$. The ordinary meaning of the above limiting statement is however that the real variable $x$ approaches {\em linearly} and {\em uniformly} with constant rate 1 to the limiting value 0 remaining always on the connected line segment $[0,x]\subset R$, $R$ being the set of real numbers. 

Let us recall also that this continuous and connected flow may, however, be broken by introducing a scaling of the form $f(x)=ax, \ 0<a<1$ and subsequently iterating  infinite number of times. The original length (Lebesgue measure)  $x$ of the closed interval $[0,x]$ now reduces to 0 successively by jumps as $a^n x, $ as $n\rightarrow \infty$. Conventionally, one feels that these are only two possibilities with unique end results. 

\subsection{Notations and Preliminary Results}

This paper is based on our previous work reported  in Ref.\cite{dp1,dp2,dp3,dp4, dp5}. In Ref.\cite{dp1, dp2, dp3}, the formalism of a scale invariant analysis is developed to give a differentiable structure on a Cantor subset of $R$. To uniformize the notations and 
for the sake of clarity, let us begin by defining and fixing our notations. By $x$ we denote a real variable. As $x\rightarrow 0$, we often suppose $x$ to define a scale. The original real variable then gets a deformed (or dressed) value $X$, and the associated infinitesimals relative to the scale $x$ are denoted by $\tilde X$ so that the inequality $0<\tilde X<x\leq X$ remains valid. The corresponding scale invariant (normalized) real variable and infinitesimals are denoted respectively by $Y=X/x$ and $\tilde Y=\tilde X/x$. The infinitesimals $\tilde X$ and the scale invariant infinitesimals $\tilde Y$ are said to live in the gaps of a fat Cantor set ${\tilde C}$ and residual renormalized Cantor set $\tilde {\cal C}$ respectively. The deformed scale invariant real variable $Y$ then resides in a measure zero Cantor set $\cal C$. In short, symbols with a tilde sign  relate to infinitesimals.

A Cantor set is a compact, perfect, totally disconnected, metrizable topological space. In this work we consider  Cantor sets those are realized as a proper subset of the real line. A Cantor set $C$  is measure zero (ie. thin) if its Lebesgue measure vanishes $m(C)=0$. Otherwise $C$ is said to be fat having a positive Lebesgue measur $m(C)>0$. A thin Cantor set $C$ is a $s$-set if the corresponding Hausdorff measure has a finite nonzero measure viz, $0<H^s(C)<\infty$. Here, $s$ denotes the Hausdorff dimension of the thin Cantor set.  A fat set $C$ is a 1-set, viz,  $H^1(C)=m(C)$. Recall that the measure of the closed interval $[0,x]\subset R$ varies as $x$. However, the Hausdorff measure of a thin $s$-set $C\subset [0,x]$ is expected to vary as $x^s$, $x\rightarrow 0^+$. Since for a fat set, Hausdorff dimension reduces to the topological dimension 1, an appropriate variant of the Hausdorff dimension is the {\em uncertainty or fatness exponent} $\beta$ (say), so that the appropriate measure is expected to behave as $x + x^{\beta}, \ 0<\beta<\infty$ (for a detailed justification from dynamical system theory approach see \cite{fermar}; we note simply here that the fatness exponent $\beta$ coincides with the uncertainty exponent $\alpha$ when $0<\beta<1$). For a more detailed justification of these facts from the point of view of the present analysis, see Ref.\cite{dp3}. We note here that for a thin $s$-set $\beta=1-s$. However, for a fat Cantor set no such identification is valid. However, as shown in \cite{dp3}, $\beta>1$ may be interpreted as the inverse of the Hausdorff dimension of a renormalized residual set $\cal C$ when the positive measure of the fat set $C$ is eliminated recursively by rescalings. The model of a measure zero Cantor set $C$ that we keep in our mind is given as the limit set of an Iterated function system (IFS) defined by $f=\{f_i| \ f_i:[0,1]\rightarrow [0,1], \ i=0,1\}$, so that $C=f(C)$, where $f_i(x)=ax+i(1-a),$ and the scale factor $a$ is defined by $2a+c=1$. Each iteration of the IFS thus removes an open interval (gap)$O_n$ of length proportional to $c$ leaving out two smaller closed subintervals $I_{nm}$ of size  proportional to $a$ each, from each of the previous level closed intervals $I_{(n-1)m}$, so that $C=\underset{n=1}{\overset{\infty}{\bigcap}}\underset{m=1}{\overset{2^n}{\bigcup}} I_{nm}$.  We also assume $a\in (0,\frac{1}{2})$, to satisfy the so called open set condition. For a fat set no IFS for a finite set of $f_i$s is usually available. However, by allowing for a variable gap of size $c_n$, a fat Cantor set is formed of measure $m(C)=\underset{0}{\overset{\infty}{\Pi}} (1-c_n)<\infty$.

A Cantor set $C$ carries a natural ultrametric (i.e., stated in otherwords, the topology of $C$ inherited from the usual topology of $R$ is equivalent to an ultrametric topology). The ultrametric structure defined in \cite{dp2, dp3} is, however, inequivalent in the sense that a measure zero set may acquire a positive measure leading to some novel applications in the context of dynamical system theory. In the following we formulate  the relevant key concepts more carefully in $R$, superseding, in our view,  those available in \cite{dp4, dp5}. 

\section{Nonlinear Analysis: Basic Concepts and Elementary Results}
We introduce nonlinearity in the flow $x\rightarrow 0$ of  $x$ by following definitions.

{\bf Definition 1:} Let $x\in I=[0,1]\subset R$ and $x$ be {\em arbitrarily small}, i.e., $x\neq 0$, but, nevertheless, $x\rightarrow 0^+$. Then there exists  $\delta>0$ and a set of (positive) \emph{relative infinitesimals} $\tilde x$ relative to the scale $\delta$ satisfying $0<\tilde x<\delta\leq x $ and the \emph{inversion} rule $\tilde x/\delta\propto
\delta/x$. The associated scale invariant (dynamic) infinitesimals are defined by $
\tilde X=\underset{\delta\rightarrow 0^+}{\lim} \tilde x/\delta$. The set of relative infinitesimals are denoted asymptotically as ${\bf 0}=\{0, \pm \delta \tilde X\}$, with $\delta\rightarrow 0$, when the scale invariant dynamic infinitesimals are denoted by ${\bf \tilde 0}= \underset{\delta\rightarrow 0}{\lim }{\{0, \pm\tilde X\}}$. When the variable $x$ itself is considered to define a scale, the relevant scale invariant quantities are denoted respectively as $Y=X/x$ and $\tilde Y=\tilde X/x$.

{\bf Definition 2:} A quantity (variable) $X$ is said to be {\em dynamical} (or dynamic) if if its {\em value} undergoes a sponataneous monotonic (ie, directed) variation in a limiting (asymptotic) problem of the form $x\rightarrow 0$ or $x\rightarrow \infty$.

As an example, for a continuous variable $x$ approaching $0^+$ and considered itself as a scale, a class of the relative infinitesimals are represented as $\tilde X \propto x^{1+l}(1+o(x)), \ 0<l<1$, corresponding to the (deformed) real variable $X=x^{1-l}$, so that the corresponding scale invariant 
infinitesimals are defined by the formula $\tilde Y=\lambda x^l +o(x^m), \ m>l$. Choosing the parameter $\lambda$ from an open set $V\subset (0,1)$, we get a class of relative infinitesimals $\tilde X$ belonging to an open subnterval of $(0, x)$, all of which are related by the inversion rule to  the real number $X^{1-l}$ for given  $x$ and $l$. Moreover, allowing $\lambda$ to vary from a countable number of disjoint open intervals $V_r, \ r=1,2,\ldots$, one can introduce a countable, disjoint class ${\tilde I}_r$ of relative infinitesimals, all of which are related by the inversion rules $\tilde X_r/\delta = \lambda_r\delta/X$, where $\lambda_r\in V_r$. Notice that $\underset{r}{\bigcup}{\tilde I}_r \subset (0,\delta)$. It  follows therefore that the infinitesimals actually reside in the gaps of an arbirarily chosen Cantor subset $\tilde C$ of the open interval $(0,x)$. As a consequence, the nonlinearity induced via the inversion rule in the apparently linear limiting motion of $x$ naturally leads to a large family of Cantor like fractal sets.  

Notice that in the limit $\delta \rightarrow 0$, the set of relative infinitesimals apparently reduces to the singleton $\{0\}$. However,  scale
invariant infinitesimals $\tilde X$ need  not be nontrivial. 
For definiteness, the ordinary zero (0) is called the \emph{stiff} zero, when
the nontrivial infinitesimals are called \emph{soft or dynamic} zeros living in the extended set ${\bf 0}/\{0\}$. The ordinary
real line $R$ is then extended over $\mathbf{R}=\{\mathbf{\tilde x}: \mathbf{
\tilde x}=x+\mathbf{0}, \ x\in R\}$, which as a field extension, and because of the Frobenius theorem, must be an infinite dimensional nonarchimedean space \cite{dp4,dp5}. 

Clearly, Definition 1 constitutes a concrete, constructive definition (realization) of infinitesimals, having, as would become evident, nontrivial influences on the very structure of real numbers (at least in the neighbourhood of 0) in the form of an ultrametric valuation.

{\bf Definition 3: Non-archimedean Norm} The mapping $v:\mathbf{0}\rightarrow I^+=[0,1]$
\begin{equation}\label{norm}
v(\tilde x) := \underset{\delta\rightarrow 0^+}{\lim}%
\log_{\delta^{-1}} {\tilde X}^{-1}, \ \tilde x=\delta \tilde X \ \in \mathbf{%
0}
\end{equation}
\noindent together with $v(0)=0$ defines  a non-archimedean (ultrametric ) norm on the set of infinitesimals $\mathbf{0}$. Infinitesimals weighted with above absolute value are called \emph{valued} infinitesimals (for a proof of ultrametricity see \cite{dp3}).

To recall, an absolute value is non-archimedean (ultrametric) if the triangle inequality is replaced by the stronger ultrametric inequality $v(a+b)\leq {\rm max}\{v(a), \ v(b)\}$. For definiteness, let us assume that the valuation $v$ is discretely valued, so that $v$ assumes values from a countable set \cite{nonarch}. The bounded set of infinitesimals is therefore associated with a family of totally disconnected Cantor sets $\tilde{ C}_p$, designated by primes $p$, those are  arranged in the form of  hierarchical branching trees (See below for more details). We also verify that for the class of infinitesimals $\tilde x \propto \delta^{1+l}(1+o(\delta)), \ 0<l<1$, introduced already following Definition 2, $v(\tilde x)=l$.

Let us now notice that the ordinary Euclidean norm for a dynamic infinitesimal would vanish trivially because $|\tilde x_r|=\delta|\tilde X_r|=0$, in the limit $\delta\rightarrow 0$. Further, there exists as such no a priori principle in the classical analysis to select or switch the motion of a variable of the form $\tilde x\rightarrow 0$ (because of the original variable $x\rightarrow 0$) to its scale invariant factor $\tilde X=\tilde x/\delta$. The nontrivial valuation now opens up a new window to initiate such a transformation.

{\bf Remark 1:} The scale invariant infinitesimals $\tilde X$ can not be a true constant. For if $\tilde X=\mu$, a constant, then the valuation of the corresponding infinitesimal $\tilde x=\mu x$ must vanish: $v(\tilde x)=\underset{\delta\rightarrow 0}{\lim} \log _{\delta}\mu =0$. As a consequence, such  a class of infinitesimals belongs to the equivalence class of the trivial infinitesimal 0. It also follows that infinitesimals are determined upto a multiple of a constant factor. Notice also that $\tilde x=\lambda \delta^{1+v(\tilde x)}$ corresponds to a scaling law on a fat Canotr set $\tilde C$ with fatness exponent $\beta=1+v(\tilde x)$, when the real variable $x$ scales as $\delta^{1-v(\tilde x)}$. Accordingly, $v$ has the status of the exterior dimension of $\tilde C$ \cite{dp3}.

{\bf Further Explanations:} As the real variable $x$ approaches 0, and becomes arbitrarily small, $x$ starts behaving as a {\em scale}, so that in relation to the scale $x$ a slightly displaced (deformed/dressed) value of the (real) variable, written in the form $X=x+h(x):=x\times x^{-v(\tilde X)} >x$, would now assume the role of the original real variable. We now denote by $\tilde X$ a class of relative infinitesimals living in a subinterval $\tilde U$ of $[0, x)$ satisfying $\tilde U\subset [0,x)\subset [0,X]$ and $X\rightarrow 0$, as $x\rightarrow 0$. We also have the inversion relation $\tilde X/x=\lambda\times (x/X) <1$, for a constant $\lambda: \ 0<\lambda<1$. By varying $\lambda$ in an open  subinterval $V$ such that $V\subset (0,1)$, we get the desired class of infinitesimals $\tilde X\in \tilde U$. Notice that $X/x=1 +h(x)/x=O(1)$ and the ultrametric norm exists provided  $v(\tilde X)=\underset{x\rightarrow 0}{\lim} \log_{x^{-1}} (x/\tilde X)\approx \underset{x\rightarrow 0}{\lim} (h(x)/x\log x^{-1})$ has a nonzero finite value of the form $v(\tilde X)=\alpha(x)<<1$, a constant with respect to the infinitesimal variable $\tilde X$: $\frac{dv}{d\tilde X} = 0$. $v$ indeed is a locally constant function in the ultrametric topology \cite{dp3}. As a consequence, $h(x)=\alpha(x)x\log x^{-1} +o(1)$ (Notice that $h(x)$, determined as above, represents a {\em nonclassical behaviour}, not available in the classical analytic representation. Classically, one expects $h(x)=O(x^2)$). For definiteness, $h(x)$ is said to denote a {\em real valued infinitesimal influence} of  {\em genuine infinitesimals} $\tilde X$ living in the ultrametric space $\bf 0$.

Now, ultrametricity of the valuation demands that $v$ will either be a constant or might behave as a Cantor function \cite{dp3,dp4} : a nondecreasing continuous function $\phi(\tilde X)$ that maps $[0,x]$ onto $[0,x]$ such that $\frac{d\phi}{d\tilde X}=0$ almost every where in $[0,x]$ (cf. Proposition 2). The set of derivative discontinuities of $\phi$ defines a Cantor set ${\tilde C}_{\alpha}$ ($\alpha$ denotes a particular realization of $v$). The infinitesimals $\tilde X$ now, by definition, live in the collection of (clopen) gaps (in the topology induced by ultrametric norm $v$) of ${\tilde C}_{\alpha}$ (cf. the example below Def.2), and the classical analytic derivative discontinuities are interpreted as sources of further variations in a neighbourhood  of $\tilde y\in \tilde {\cal C}$, where $\tilde y=\tilde X/x$ denotes a scale invariant variable in [0,1] and $\tilde{\cal C}$ is the corresponding image Cantor set under the mapping generated by $\tilde y$.

As stated above,  we established that $v$ can be identified with a locally constant Cantor function with double logarithmic variability \cite{dp3, dp4}. Improving this result we now have following two propositions.

{\bf Proposition 1:} The derivative discontinuity of a Cantor function is removed by smooth jumps revealing a nonclassical (double logarithmic) variability in the valuation in the neighbourhood of point of a Cantor set.

{\em Proof}: Let $\tilde y_0\in \tilde{ \cal C}_{\alpha}$. Because of total disconnetion and the density of $\tilde{ \cal C}_{\alpha}$ in [0,1], there exist arbitrarily small clopen gaps ${\cal J}_-$ and ${\cal J}_+$ respectively in the left and right hand sides of $\tilde y_0$. Let $\tilde y_-\in {\cal J}_-$ and $\tilde y_+\in {\cal J}_+$. By definition, $\tilde y_{\pm}= \tilde X_{\pm}/x=x^{\mp v(\tilde X_{\pm})}:={\tilde y_0}^{\mp v_{\pm}}$, where $\tilde X_{\pm}$ lie in the appropriate gaps $J_{\pm}$ of ${\tilde C}_{\alpha}$. Further, classical analytic obstruction to a smooth variability in the neighbourhood of $\tilde y_0$ is now removed by invoking transition via smooth jumps of the form $\tilde y_+={\tilde y_-}^{-c}$ for a $c>0$, so that we have $c=v_+/v_-$. But, for an infinitesimally small change, we may set $v_+=v_-\times x^{-\psi(x)}$ so that the smooth variability of $v(\tilde X)$ around $\tilde y_0$, say, is controlled by the scale invariant equation 
\begin{equation}\label{sie}
\log x^{-1}{\frac{d\psi}{d\log x^{-1}}}= \psi
\end{equation}
(so that $c$ is essentially is locally constant viz, $\frac{d c}{d x}=0$) where $x$ now denotes a scale invariant variable that lives and varies by smooth jumps respectively, in and out of a gap, from  a countable collection of gaps of another, arbitrarily chosen Cantor set $\tilde{\cal C}_{\alpha^{\prime}}$ that replace the given Cantor  set point $\tilde y_0$. As a consequence, a (finite) jump from $\tilde y_-$ to $\tilde y_+$, as designated above, actually is reduced to infinitely many infinitesimal jumps,  thus revealing , in turn, a new level of nonclassical variability in the said neighbourhood. $\Box$

{\bf Proposition 2:} The ultrametric valuation $v(\tilde X)$ has the general form $v(\tilde X)=ax +b x^{\beta}, \ 0\leq a,\ |b|<1$, $ \beta>1$ being a fatness exponent, ie, $v$ is either  a constant or a locally constant Cantor function, or a combnation of both as above, with respect to $\tilde X$. The variability of $v$, however,  is exposed in a double logarithmic scale. Further, $v$ represents a multifractal measure. 

{\em Proof:} By construction, 0 of $R$ is replaced by a Cantor subset $\tilde C_{\alpha}$ of $[0,x]$ when $x\rightarrow 0^+$. Let $\tilde X\in \tilde U$, an open gap, in the usual topology. We set $v(\tilde X)=ax$, when $\tilde U$ belongs to the subclass of countable gaps that are deleted by an underlying IFS (finite or infinite) leading to the Cantor set ${\tilde C}_{\alpha}$. This {\em unique} choice is guided by the fact the open interval $\tilde U\subset [0,x]$ is of Lebesgue measure $ax<x$ for a suitable $a$. Further, it also tells that $\tilde C_{\alpha}$ is a positive measure (fat) Cantor set with measure $ax$ and fatness exponent  $\beta$, say. When, on the other hand, $\tilde X$ is drawn from  an open interval $\tilde U_{\infty}=(a_{\infty},b_{\infty})$ whose end points $a_{\infty}$ and $b_{\infty}$ are realized as limit points of sequences of end points of the countable family, we set $v(\tilde X)=b(a) x^{\beta}$, $\beta=1+s_0(x)$ being the fatness exponent for the fat Cantor set ${\tilde C}_{\alpha}$. The  reason for this choice corresponds to the fact that $v$ now represents a Cantor function (for a proof \cite{dp3,dp4}) with the desired functional dependence (Proposition 1) when $s_0=s_0(a,x)$ denotes a variable Hausdorff dimension of the renormalized residual Cantor set $\tilde{\cal C}_{\rm res}$ when the positive measure of underlying fat Cantor set ${\tilde C}_{\alpha}$ is removed (c.f. Sec.1.1). More generally, when $\tilde X$ is drawn from either of an open gap $\tilde U$ or of a limiting gap of the form ${\tilde U}_{\infty}$ (this indistinguishablity arising from the inherent limiting motion $x\rightarrow 0$) $v(\tilde X)$ can only be written as a linear combination of above two cases, when the coefficient $b$ is allowed to vary in (-1,1) with the restriction $v(\tilde X)>0$. Since, an open interval of points are assigned a constant (relative to the infinitesimal variable $\tilde X$) value, $v$ is  trivially an ultrametric over each  gap from the family of gaps designated as $[0,x]/\tilde C_{\alpha}$. However, recalling $(0\approx) \ x^2\propto X\tilde X={\tilde X}^{\gamma(\tilde X)}$, a constant, both $\tilde X$ and $X$ enjoy small scale variability that gets exposed, thanks again to Proposition 1,  in double logarithmic scale $\log\log\tilde X$. Rewriting $v$ using $\tilde X$, the doulbe logarithmic variation follows. Finally, the Cantor set  $\tilde C_{\alpha}$ and the related valuation $v$ both enjoy multi- parameter arbitrariness, thus establishing that $v$ defines a multifractal measure over a class of Cantor sets $\tilde C_{\alpha}$. $\Box$

{\bf Corollary 1:} {\em \cite{dp2,dp3}}  Let $x\rightarrow 0^+$ on a  measure zero Cantor set $C$. Then the valuation for the associated infinitesmals has the form $v(\tilde x)=b x^s$, $s\in (0,1)$ being the Hausdorff dimension of $C$ and $b$ being a (true) constant.

{\bf Corollary 2:} The valuation of Proposition 2 can be written more poignantly in the form 

\begin{equation}\label{meas}
v(\tilde X_p)= a_p x + b_{p}x\times Y
\end{equation}
where $Y= (\frac{\bar x}{x})^{s}\geq 1$ and $\bar x$ is a real variable that lies in the infinitesimal line segment that replaces the point at $x$, and $s$ is interpreted as the Hausdorff dimension of the (measure zero) residual Cantor set $\tilde {\cal C}_{p}$ over which the fluctuating motion takes place once the linear connected flow is discarded recursively by rescalings (we henceforth denote $\tilde{\cal C}_{\rm res}$ by $\tilde {\cal C}_p$, and the symbol $p$ replaces the suffix $\alpha$).

{\bf Remark 2:} The above typical behaviour of a multifractal measure is recently observed and studied in quantum chaos \cite{weyl} (see also \cite{fermar}). The proof of Cor.2 follows  from the preliminary remarks in Sec.1.1 on the fatness index $\beta$ and also exposes how a linear connected line segment might arise dynamically in the vicinity of a point of a fat Cantor set. Indeed, for a $\beta$ such that $1<\beta<2$, we set $\beta=1+s, \ 0<s<1$ so that the measure of prop. 2 behaves  as
$v(\tilde X_p)= a_p x + b_{p}x\times (\frac{\bar x}{ x})^{s}$. As a consequence, the original motion of $x$ is transferred to $v(\tilde X(x))$ consisting of a linear motion over a connected line segment, along with a local fluctuating motion on a residual Cantor set.

{\bf Corollary 3:} The ultrametric valuation  has the more general form  $v(\tilde X)= ax(1+\sum b_i x^{s_i}), \ 0<a, \ |b_i|<1$ and the sequence of exponents $0<s_{i+1}<s_i<1$ corresponds to finer scales residual Cantor sets those are reached successively by rescalings as $x\rightarrow 0^+$.

The proof of ths corollary becomes obvious in the light of Proposition 2 and Theorem 1 (Sec.3).

Before leaving this section, let us remark on the metric structure of the extended set $\bf R$.

{\bf Definition 3:} The valuation $v$ on $\bf 0$ induces an ultrametric norm $||\cdot||$ on the whole set $\bf R$ with the property $||X||=v(X)$, when $X\in \bf 0 $, but $||x||=|x|$, the usual Euclidean norm, when $ x\in \ R-\{0\}$.

{\bf Lemma 1:} $||\cdot||$ is an ultrametric norm on $\bf R$.

{\em Proof:} Ultrametricity of $||\cdot||$ follows from $v$, and the observation that for $X\in {\bf R}-\{0\}$, $||X||=||x+\tilde x||, \ \tilde x\in \bf 0$, so that $||X||=\max \{||x||,v(\tilde x)\}=|x|$. $\Box$

\section{Relation with Arithmatical functions}

We now study the variability of the multifractal measure $v(\tilde X)$ in more details. Let us recall that the dressed (deformed) real variable, written scale invariantly, has the form $Y=X/x=x^{-v(\tilde X)}$. Classically, one expects $Y=1$ identically for all $x$. In the present formalism, $Y$, however, experiences a variation, as the original variable $x$ continues to probe the space of infnitesimals deeper and deeper in an hierarchical manner.

{\bf Theorem 1:} The set of infinitesimals $\bf 0$ equipped with the valuation $v$ is an infinite dimensional ultrametric space. Equivalently, it is a rooted, directed tree.

{\em Proof:} 
Infinitesimals, represented as ${\bf 0}\ni\tilde X_p\propto x\times x^{v(\tilde X_p)} +o(x)$, as $x\rightarrow 0$, using a scale $0\approx x\in R$, are {\em new } elements in $\bf 0$, which can not be accommodated in the ordinary real line $R$. In fact $\bf 0$, when restricted to $R$, reduces to the singleton $\{0\}$. However, the scale invariant components $\tilde Y_p=\tilde X_p/x= x^{v(\tilde X_p)} +o(1)$, are nevertheless, nontrivial, and so the uncountable set ${\bf 0}/\{0\}$ must be raised to a higher dimensional (normed) vector space. Indeed, $\bf 0$ is a vector space under usual vector addition and scalar multiplication over the field of real numbers $R$.  

Next, the open set of infinitesimals as a subset of the open interval $(0,x)$ for each fixed $x$, is realized as the union of countable number of disjoint open subintervals $\underset{p}\bigcup U_p\subseteq \bf 0$. To each open set $U_p$, associate a class of infinitesimals $\tilde X_p$, weighted with valuations, for instance, $v(\tilde X_p)=a_p x$ (by Prop. 2). Further, $\tilde X_p\in U_p \ \Rightarrow \ \mu\tilde X_p\in U_p$ ($\mu$ beng a true constant), since $v(\mu\tilde X_p)=v(\tilde X_p)$. As a consequence, $\bf 0$ is a normed vector space over $R$ with trivial ultrametric norm $|\mu|=1$ for $\mu\in R$.

To justify this more precisely, let $e_p=(0,0,\ldots 1, 0\ldots)$, with 1 at the $p$-th place, denote the standard basis for an infinite dimensional space. Let $\tilde {\bf X}_p=\tilde X_p\cdot e_p$ denote an (countably) infinitesimal vector along $p$-th direction. Hence, any vector $\tilde X\in \bf 0$ has the representation $\tilde {\bf X}=\sum \tilde X_p\cdot e_p$. Define a norm $||\cdot||_0$ on $\bf 0$ equipped with the valuation $v$ by $||\tilde {\bf X}_p||_0=v(\tilde X_p)|| e_p||$ with the condition $||e_p||=1 \ \forall p$. Then $||\tilde {\bf X}||_0= \sup v(\tilde X_p)$, because of the ultrametricity. We also verify that the ultrametric convergence condition $\lim v(\tilde X_p)=0$ as $p\rightarrow \infty$ holds good, so that the infinite summation for $\tilde {\bf X}$ is well defined. As a consequence, $\bf 0$ is an infinite dimensional ultrametric vector space.  

This set of infinitesimals would align as an orthogonal subspace ${\bf 0}_{v}$, relative to the horizontal real axis, with a vertex point (node) at the origin 0. Since, ${\bf 0}$ is an ultrametric space, and since, every ultrametric space is equivalent to a directed tree, ${\bf 0}$, is a directed tree with a root. $\Box$

{\bf Lemma 2:} The indexing parameter $p$ varies over the set of primes, and the scale factors appearing in the explicit formula for $v$ in Prop. 2 are inverse powers of prime $p$.

{\em Proof:} Let $\delta=x_0$ be a scale, where $x_0$ denote a rational approximant of $x$. Clearly, $v(\tilde x)$ would vanish, if one restricts the motion of $x$, as it approaches 0, only on the real number system $R$. For nontrivial values, $x$, and hence $\delta$, must flow through a neighbourhood of 0, consisting of {\em nonreal} elements. This is realized by an (infinite dimensional) ultrametric extension of $R$ on to $\bf R$, in which ordinary singelton set $\{0\}$ is replaced by the fattend space ${\bf 0}$ which is defined symbolically as  $\log {\bf 0}/\{0\}= \Pi Z_p$, where $Z_p$ is the ring of $p$-adic integers, so that the logarithm of a scale invariant infinitesimal $\log \tilde {\bf X}\in \Pi Z_p$ (c.f. \cite{dp4}). By Ostrowski's theorem, the ultrametric norm on $\bf 0$, when restricted over each of $Z_P$'s, would reduce uniquely to the $p-$adic norm $|\cdot|_p$. As a consequence, $v(\tilde X_p)\approx \alpha_px$ (cf. Prop. 2), when $x$ takes the first inversion induced jump at the scale $\delta$ and lands on the subinterval $U_p\subset (0,x)$, where the scale factor $\alpha_p$ is a rational of the form $\alpha=p^{-n}$ (in the limit $x\rightarrow 0$, $\alpha_p$ is raised to a $p$-adic number $\hat{\alpha_p}$ such that $|\hat{\alpha_p}|_p=p^{-n}$). Let $\xi: Z_p\rightarrow {\cal C}_{p}$ denote the natural homeomorphism between $Z_p$ and a Cantor subset ${\cal C}_{p}$ of [0,1] so that $\xi(\hat{\alpha_p})=\alpha_p$ and $||\hat{\alpha_p}x||=v(\hat{\alpha_p})|x|=p^{-n}x$. The logarithm of the corresponding deformed scale invariant variable $Y=X/x$ thus belongs to the family of Cantor sets ${\cal C}_{p}$ (i.e., $\log Y\in {\cal C}_{p}$), which are accessed one after one by rescalings $x\rightarrow Y=X/x$, as $x\rightarrow 0$ more accurately through the $p-$adically generated scales $x_p=p^{-n}x, \ p\rightarrow \infty$. $\Box$

To summarise, in a limiting problem with $x\rightarrow 0$,  infinitesimals $\tilde X_p$ living  in ${\bf 0}_{v}$ are activated, so as to infuse a fattening influence on the ordinary real analytic variables of the form $x$, which gets deformed as $X$, so that the ordinary scale invariance of the equality $Y=X/x$ is violated, giving rise to a new nonarchimedean scale invariance of the form $Y=X/x=x^{-v(\tilde X)}, \ \tilde X\in \ {\bf 0}_{v}$. We next evaluate the actual value of  $v(\tilde X)$.

{\bf Theorem 2:} $v(\tilde{\bf X}) = x \times \Pi(x^{-1}) $, where $\Pi(x^{-1})$ is the prime counting function \cite{rh}.

{\em Proof:}  As $x\rightarrow 0$ and gets arbitrarily small, the linear motion of $x$ gets frozen, leading to a scale $x$, and thereby making room for nonlinear hopping motions for the deformed variable $X\geq x$, as represented by the deformation formula $Y=X/x=x^{-v(\tilde X)}$ over the family of residual Cantor sets ${\cal C}_{p}$. Accordingly, the original motion is transferred to the hidden ultrametric space of infinitesimals $\tilde X$; an infinitesimal now continues to vary over several branches of the infinite adelic tree of the form $\Pi Z_p$. Each jump from one $Z_p$ to the next $Z_q$, $q$ being the immediate successor of $p$, would induce a mirror variation (by induced homeomorphisms)  from ${\cal C}_{p}$ to ${\cal C}_{q}$ (via the map $\xi$). The valued measure $v(\tilde X)$ will thus have the form, after having initiated the first level transition from linear motion of $x$ on $R$ to  $\tilde {\cal C}_{2}$ (recall that $\tilde {\cal C}_{2}$ belongs to the orthogonal subspace of the homeomorphic image of ${\bf 0}_v$), given by $v(\tilde X)= 2^{-n}x(1 + y_2), $ where $y_2=b_2x^{-s_2} \ \in \tilde {\cal C}_{2}\subset [0,1]$. Let, $x_1=2^{-n}x$ and choose, $b_2=2^{(-n-r)s_2}$, so that $v(\tilde X)= x_1(1 + y_2), \ y_2=2^{-rs}x_1^{-s_2}$, and $y_2<<1$, initially though,  is a slowly growing real variable lying in the connected line segment (0,1) (by Remark 2) obtained by the defining rescaling (of $y_2$ ) from a sufficiently small connected neighbourhood of an appropriate point of $\tilde {\cal C}_{2}$. So, when we take the limit $x\rightarrow 0$, the first order inverted rescaled variable $y_2$ would grow to a value $\lesssim O(1)$ by linear translation and then undergo the next level transition to the branch $\tilde {\cal C}_{3}$ by inversion $(1\gtrsim ) y_2\rightarrow y_2^{-1}=1 + y_3$ where the second order rescaled variable $y_3$ now lives in and subsequently grows to O(1) in $\tilde {\cal C}_{3}$, so as to make yet another transition to the next branch $\tilde {\cal C}_{5}$ following the rule $(1\gtrsim ) y_3\rightarrow y_3^{-1}=1 + y_4$ with $y_4$, initially a small rescaled variable would again grow slowly to O(1) to initiate yet another inversion induced transition to $p$-adically generated residual Cantor set $\tilde {\cal C}_{7}$ and so on, yielding the following universal pattern of (internal) variation in the valuation of the form $ v(\tilde X)=x\times (1 + 1 + \ldots y_n)$ (we replace $x_1$ by $x$ and understood it as a rescaled variable living in (0,1) and approaching $0^+$, because of the inherent scale invariance) where $y_n$ belongs and varies on the $n$th $p$-adic induced Cantor set $\tilde {\cal C}_{p}$. This proves the theorem once we note that $\Pi(x^{-1})=\underset{p<x^{-1}} {\sum} 1$.  $\Box$

One also deduces an interesting consequence of the above result.

{\bf Proposition 3:} The classical value of $\log X/x=0$ gets a non-classical counterpart $\log X/x=1 + O(x^{\frac{1}{2}-\epsilon})$ for $\epsilon>0$.

{\em Proof:}  The deformation caused by nontrivial valuation is defined by $\log X/x= v(\tilde X) \log x^{-1}$. Theorem 2 tells that the limiting value of this equality of course should be 1 viz, $\underset {x\rightarrow 0}{\lim}\log X/x=1$, thus proving the prime number theorem $\underset {x\rightarrow 0}{\lim} x\log x^{-1}\Pi(x^{-1})=1$. Putting $y=x^{-1}$ we get the prime number theorem in the conventional form $\underset {y\rightarrow \infty}{\lim} (\frac{\log y}{y}) \Pi(y)=1$. The asymptotic form of this limiting equality now has the form $(\frac{\log y}{y}) \Pi(y)=1 +O(y^{-\alpha}), \ 0<\alpha<1$.

To justify the above equality and also determine the exponent $\alpha$ more accurately, let us restate the above again using the variable $x$: $x\log x^{-1}\Pi(x^{-1})=1 +O(x^{\alpha})$. We first notice that the rhs is just another way of writing the non-classical deformation $\log X/x= 1 + O(x^{\alpha}) $. Further, the inversion rule $\tilde X X\propto x^2, \ x\rightarrow 0^+$ tells $X\sim x^{\frac{1}{2}-\epsilon_0}$ (so that $\tilde X\sim x^{\frac{3}{2}}+\epsilon$ and $x^{\epsilon_1-\epsilon_0}\approx 1$), for an $\epsilon_0>0$ ($\epsilon_1>0$), leading to the most natural estimate $\alpha=\frac{1}{2}-\epsilon$ for an $\epsilon >0$.   $\Box$

{\bf Remark 3:} We interprete the above paradoxical result as follows. Let us recall that vacuum in classical mechanics is a mere emptiness, while that in quantum mechanics is a sea of elementary particle-antiparticle ( for instance, electron-positron) pairs. According to Dirac, electrons do have stable physically realizable states when positrons mostly arise as ``holes" in the electron-positron sea. Recall also that physical states of a positron are highly unstable. Exactly in a similar vein, the present nonlinear, nonclassical analysis now offers a bigger panorama to view the number thoeretic vacuum, viz, the number zero (realized here as the limit of the closed interval $[0,x],$ as $x\rightarrow 0$) or the empty set (ie, the limit of the open interval $(0,x)$). The ordinary number 0 or the sense of emptiness is now realized as an infinite dimensional ultrametric space of  infinitesimals accommodating genuine  infinitesimals $\tilde X$ (analogue of positron holes) and their real valued images by inversions, the so called real valued infinitesimals (equivalent to the physical states of an electron), denoted here as $X$, so that the scale invariant balancing equality $\tilde X X\propto x^2$ holds good dynamically in the limit $x\rightarrow 0$. In short, the framework of the present nonlinear analysis can support emergence of complex structures quite naturally in contradistinction with the classical analysis.

\section{Alternative Justifications and Naturalness}

Classically, the statement $\underset{x\rightarrow 0}{\lim} x=0$ means that both the Lebesgue measure and cardinality of the open interval $(0,x)$ vanish in the limit $x\rightarrow 0$. The cardinality thus has a huge singularity: an uncountable set is reduced suddenly to the empty set. Classical analysis fails to answer the question regarding  whereabouts of the uncountable number of real points (numbers) in the interval $(0,x)$ as $x\rightarrow 0$.  The above non-classical results and insights now offer a resolution of these paradoxical properties. This also provides, a first principle explanation, even on the basis of the classical analytic framework, of the growth of values and measures as discussed  in the derivation of the prime number theorem. The naturalness of the definition of the ultrametric norm (Definition 3) is also demonstrated afresh.

Let us consider the scale invarinat equality (a conservation principle)

\begin{equation}\label{consv}
x=\{x a(x)^{n\log x^{-1}}\}\times p^{n\log X}, \ 0<x<1, \  X>1
\end{equation}
which clearly formulate the precise classical analytic setting how a linear flow, as it were, could be  transformed into a scale invariant nonlinear flow in the logarithmic scales. Notice that as $x$ approaches 0 linearly (the lhs of the equality), a new level of reduction of the linear measure 1 is initiate at the rhs, by a scaling factor $a(x)$ depending on $x$. As a consequence, the reduction of the linear measure of scale invariant value 1 (on the lhs) at each scale by the factor $a(x)^{n\log x^{-1} }$ is corrected multiplicatively by creating an {\em extra space} in a vertical direction containing the fattened variable $X>1$ in the logarithmic scale for another scaling variable $p>a^{-1}$ and as a consequence, in the limit $n\rightarrow \infty$, and hence as $x\rightarrow 0$, the linear variable $x$ is transformed into a nonlinear flow characterised by the scaling law $X^{-1}=x^{s(x)}$, where $s(x)=\frac{\log a(x)^{-1}}{\log p}<1$ is the (local) fractal (Hausdroff) dimension of a Cantor like set that might have  been assumed to form spontaneously in an extended neighbourhood of 0, because of the limit process. We note that, because of the inversion rule, the original variable $x$ on its approach towards 0, would get reversed (ie, experience a bounce) close to 0, so that the corresponding  infinitely large fattened variable $X_s= \tilde X^{-1}= xX^{-1} $ will live and vary on a fat  Cantor set $\tilde C$ of fractal dimension, (actually, the fatness exponent) $\tilde s=1 +s(x)$, when $s(x)$ may now be interpreted as the exterior dimension of $\tilde C$. This, in turn, would have a very important consequence: {\em the value 0 is actually never reached; the linear flow along the positive $x$ axis, in the interval $(0,x), \ x\downarrow$, will be transferred in a scale invariant manner, to a new branch in the inverted variable $X_s\sim x^{\tilde s}$, by inversion induced jumps. The subsequent (linear) flow of $X_s$ toward 0 will then again experience a bounce, and so on leading to a very intricate multifractal tree like structure.} The multifractality follows from the existence of a one parameter family of multiplicative scale factors $p(a)>a(x)^{-1}$ depending on $a(x)$ and hence on $x$.

The above arguments leading to the infinitesimal (higher dimensional) space that is revealed on the logarithmic scales, also offers us {\em a mechanism to reconstruct} the definition of the scale invariant valuation in a rather novel manner, thus establishing naturalness of Definition 3. 
 
 To this end, we rewrite the conservation equation (\ref{consv}) of measure in an alternative manner so as to get the slightly modified version in the form

\begin{equation}\label{consv2}
x=\{xe^{-nv(\tilde X)\log x^{-1}}\}\times e^{n \log  Y},
\end{equation}
To restate once more, as the real variable $x$ approaches 0 continuously (lhs of the above equation), its motion is momentarily frozen to $x$ ( say) (rhs of the equality)  at any instant of time $t$ (say), thus  making a room for a new level of variability by a nontrivial scaling with logarithmic exponent $v(\tilde X)\log x^{-1}$, which must have to be balanced by bringing in another scaling variable $p>1$, for an exponent corresponding to a new scale invariant variable $Y>1$ living in an orthogonal subspace of the extended neighbourhood of 0. Notice that $v$ here carries the signature of invisible infinitesimals, so that $v(\tilde X)=0$ implies that $ Y=1$. For a nontrivial valuation, we would indeed recover the original definition of $v$, viz, $v(\tilde X)= \log_{x^{-1}} Y$, in the limit as the scale $x\rightarrow 0$. To make contact with  Definition 3, we note that when $x$ is identified with a scale, we should consider a slightly deformed (displaced) value of the form $X=x+h(x)$, where $h(x)$ $\sim$ O($x^{\beta})$, $ 1<\beta\leq 2$, (in ordinary analysis, of course, $\beta=2$ uniquely) so that $Y=X/x$. As a consequence, in this extended scale invariant framework, $h$ appoximately has the form $h(x)\approx v(\tilde X)x\log x^{-1}$ (in perfect agreement with what we had already established.)

As an application of the above novel classical analysis, leading, nevertheless, to nonclassical results, let us point out a realization of Chebyshev's $\psi(x)$ function $\psi(x)=\underset{p^n<x^{-1}}{\sum} \log p$ \cite{rh}. Let us work with scale invariant  real and infinitesimals defined by $0<\tilde Y<1<Y$ where $X=xY$ is the deformed variable relative to the scale $x$. $Y$ approaches 0 through scales $\delta_{np}=p^{-n}$. We set the uniform value $v(\tilde X_{np})=x/n$ for each fixed $p$. Then (\ref{consv2}) (indeed, an extended form) is rewritten as 
\begin{equation}\label{psi}
1=\{e^{-m\underset{p^n<{x^{-1}}}{\sum}v(\tilde X_{np}){\frac{\log \delta_{np}^{-1}}{\log x^{-1}}}}\}\times e^{m \log_{x^{-1}} Y}
\end{equation}
so that $Y=x^{{\frac{x}{\log x^{-1}}}\psi(x^{-1})}$ which compares correctly with Theorem 2, in the limit $x\rightarrow 0$, when we recall $\psi(x)\sim x$ for $x\rightarrow \infty$.


\begin{thebibliography}{99}
\bibitem{dp1} S Raut and D P Datta, Analysis on a fractal set, {\em Fractals}, {\bf 17}, (2009), 45-52, erratum ibid, {\bf 17}, (2009), 547.
\bibitem{dp2} S Raut and D P Datta, Nonarchimedean scale invariance and Cantor sets, {\em Fractals}, {\bf 18}, (2010), 111-117.
\bibitem{dp3} D P Datta, S. Raut  and A Raychaudhuri, Ultrametric Cantor sets and growth of measure, {\em p-Adic Numbers, Ultrametric Analysis and Application}, {\bf 3}, (2011), 7-22.
\bibitem{dp4} D P Datta  and A Raychaudhuri, Scale invariant analysis and the prime number theorem, {\em Fractals}, {\bf 18}, (2010), 171-184.
\bibitem{dp5} D P Datta, On a new proof of the Prime Number Theorem, ArXive:1103.6176v1[Math.Gen] (2011). 
\bibitem{nonarch} P Schneider, Nonarchimedean Functional Analysis, Springer, (2001).
\bibitem{fermar}J D Farmer, Sensitive dependence on parameters in nonlinear dynamics, {\em Phys. Rev. Lett}. {\bf 55}, (1985), 351-354.
\bibitem{weyl} M. E. Spina, I. Gracia-Mata, M. Saraceno, Weyl law for fat fractals, {\em J. Phys. A: Math. Theor.},  {\bf 43}, (2010), 392003.
\bibitem{rh} H. M. Edwards, {\em Riemann's Zeta Function}, Dover, (1974).
\end{thebibliography}
\end{document}